# Sharp upper bounds on the number of the scattering poles

Plamen Stefanov*
Department of Mathematics
Purdue University
West Lafayette, IN 47907

## 1 Introduction and Main Results

Let $P$ be a compactly supported perturbation of the Laplacian in $\mathbf{R}^n$, $n$ odd, defined by the "black box scattering" formalism, i.e., $P = -\Delta$ outside $B(0, R_0)$ and $P$ satisfies the hypotheses in [SjZ2], see section 2. As usual, $P^\sharp$ denotes a reference operator on the "perturbed torus" $\mathbf{T}_R$, $R > R_0$, see next section. Let $N(r)$ be the number of resonances (scattering poles) of $P$ with modulus less than $r$. One of the basic questions in the theory of resonances is to estimate $N(r)$ and, if possible, to find an asymptotic formula, as $r \to \infty$. In a pioneering work, Melrose [M1] showed that $N(r)$, related to $P = -\Delta + V(x)$, where the potential $V$ is compactly supported, has at most polynomial growth, and in an unpublished note later he improved this to $N(r) \leq Cr^{n+1}, r > 1$. Then he showed [M2] that

$$N(r) \leq Ar^n \tag{1}$$

in obstacle scattering. M. Zworski [Z2] proved (1) for compactly supported potentials. The case of elliptic second order $P$'s was resolved by G. Vodev in [V1], and in [V2] for non-self-adjoint operators. In a general black-box setting, a generalization of (1) was proved by Sjöstrand and Zworski [SjZ1]. Similar bounds are known in the semiclassical case, see e.g., [PZ2] for references. Bounds on a modified version of $N(r)$ in even dimensions were studied in [I], [V3], [V4].

It is known that the distribution of the scattering poles in various neighborhoods of the real axis depends on the geometry of the scatterer, respectively on the properties of the Hamiltonian flow associated with $P$. We will not give full account of those results and will mention only [SjZ2], [Ze] where scattering poles in sectors $0 < -\arg \lambda \ll 1$ are studied, and [S2] for upper bounds in $\{|\lambda| > 1; \Im\lambda \leq |\lambda|^{-N}\}$, $N \gg 1$.

At present, very little is known about a possible asymptotic formula for $N(r)$. In the 1D case, for $P = -\mathrm{d}^2/\mathrm{d}x^2 + V(x)$, it is known [Z1] that

$$N(r) = \frac{2a}{\pi}r + o(r), \tag{2}$$

where $a$ is the diameter of the support of $V$. For $n \geq 3$ odd, and $P = -\Delta + V(x)$, M. Zworski [Z3] proved that

$$N(r) = K_n R^n r^n (1 + o(1)),$$

under the assumption that $V \in C^2$ is radial, supported in $B(0, R)$, and $V$ has a jump at $|x| = R$ (see also Theorem 3 below). His proof also implies an asymptotic of the same kind with a different constant for the sphere resonances with Dirichlet boundary conditions (see also Theorem 1 below). The constant $K_n$ in (2) is not specified in [Z3][1].

The purpose of this work is to find an explicit constant $A_n$ such that

$$N(r) \leq A_n r^n + o(r^n) \tag{3}$$

---

*Supported by NSF Grant DMS-0400869

[1] Actually, an integral representation of $K_n$ of the kind we obtain is implicit in [Z3], see the proof of Theorem 3.

for various $P$'s, and to show that in some special cases, including those above, $A_n$ is sharp because then (3) turns into an asymptotic. The term $A_n r^n$ would then serve as a candidate for the leading term in the asymptotics of $N(r)$, if the latter exists, at least if the scatterer is spherically symmetric. In the case studied in [Z3], one can see that the constant $K_n R^n$ depends on the size of the support of $V$ only, and not on $V$ itself. This corresponds well to the known fact that the scattering determinant $s(\lambda)$ related to general $P$'s, admits an estimate of the kind $|s(\lambda)| \leq C'_n \exp(A'_n |\lambda|^n)$ in the "physical plane" $\Im \lambda > 0$ [PZ2].

Since we use Jensen's type of equality, this forces us to work with a regularized version $M(r)$ of $N(r)$, instead of $N(r)$:

$$M(r) = n \int_0^r \frac{N(t)}{t} \, dt = n \sum_{|\lambda_j| < r} \log \frac{r}{|\lambda_j|}, \tag{4}$$

where $\lambda_j$ are the resonances. The factor $n$ above can be explained by the following: $M(r)$ has asymptotic if and only if $N(r)$ has asymptotics, and then the leading terms coincide, see Lemma 1. In all cases, $N(r) \leq eM(r)$ but then the factor $e$ probably makes the estimate for $N(r)$ non-sharp.

In order to state our main results, introduce the function

$$\rho(z) = \log \frac{1 + \sqrt{1-z^2}}{z} - \sqrt{1-z^2}, \quad |\arg z| < \pi,$$

see section 4 for more details. We denote $[-\Re\rho]_+ = \max(-\Re\rho, 0)$. Then in $\mathbf{C}_+ := \{\Im z > 0\}$, the function $[-\Re\rho]_+$ is supported outside an eye-like domain $\mathbf{K}$, see Figure 1 in section 4.

We study first the case of $P = -\Delta$ in $|x| > R_0$ with Dirichlet boundary conditions. The resonances in this case are well known to be the zeros of $H^{(1)}_{l+n/2-1}(\lambda R_0)$, $l = 0, 1, \ldots$ with multiplicities equal to the dimension of the corresponding spherical harmonics eigenspace and their asymptotics follow from Olver's uniform asymptotics of Hankel's functions, see Figure 2 in section 7. The asymptotic of the counting function however, to our best knowledge, has not been studied except for [Z3], as mentioned above.

**Theorem 1** *Let $N_{R_0 S^{n-1}}(r)$, $n$ odd, be the counting function of the resonances for the exterior Dirichlet problem for the sphere $R_0 S^{n-1}$, $R_0 > 0$. Then*

$$N_{R_0 S^{n-1}}(r) = A_{S^{n-1}} R_0^n r^n + o(r^n), \quad \text{as } r \to \infty,$$

*where*

$$2 \frac{\operatorname{vol}^2(B(0,1))}{(2\pi)^n} + A_{S^{n-1}} = \frac{2n}{\pi(n-2)!} \int_{\Im z > 0} \frac{[-\Re\rho]_+(z)}{|z|^{n+2}} \, dx \, dy, \quad z = x + iy. \tag{5}$$

*In particular, if $n = 3$,*

$$\frac{4}{9\pi} + A_{S^2} = \frac{6}{\pi} \int_{\Im z > 0} \frac{[-\Re\rho]_+(z)}{|z|^5} \, dx \, dy.$$

The same results holds for Neumann or Robin boundary conditions, as it can be seen from the proof.

Numerical experiments based both on direct counting all resonances with modulus less that 67 ($N_{S^2}(67) = 522,772$), and on a numerical computation of the integral above, show that $A_{S^2}$ is in the range $(1.73, 1.75)$. Another integral representation of $A_{S^{n-1}}$ is given in Lemma 4 below.

We find it more convenient to work with a reference operator $P^\sharp$ equal to $P$ in $B(0, R)$, where $R > R_0$ is fixed, with Dirichlet boundary conditions on $|x| = R$ (to be more precise, "in $B(0, R)$" means on $\mathcal{H}_{R_0} \oplus L^2(B(0, R) \setminus B(0, R_0))$), see next section). In most interesting cases, taking the limit $R \to R_0$ would provide the best estimates. From now on, $P^\sharp$ is that operator, and $R > R_0$ is fixed.

Recall that one of our assumptions is that $N^\sharp(r+1) - N^\sharp(r) = o(r^n)$, see section 2. One can impose this assumption either on the torus reference operator, or on the ball reference operator.

For convenience, set

$$\tau_n = (2\pi)^{-n} \operatorname{vol}^2(B(0,1)). \tag{6}$$

Next theorem is the main result of the paper. It gives an upper bound for $M(r)$ for general $P$'s.

**Theorem 2** *Let P satisfy the black-box assumptions in the ball $B(0, R_0)$ described in section 2, n odd, with a reference operator $P^\sharp$ in $B(0, R)$, with arbitrary but fixed $R > R_0$. Then*

$$\left|M(r) - 2\left(N^\sharp(r) - \tau_n R^n r^n\right)\right| \leq (2\tau_n + A_{S^{n-1}}) R_0^n r^n + o(r^n), \quad \text{as } r \to \infty. \tag{7}$$

*In particular,*

$$M(r) \leq 2N^\sharp(r) + R_0^n A_{S^{n-1}} r^n + o(r^n), \quad \text{as } r \to \infty.$$

One can interpret the result above as follows: the number of resonances is bounded by the number of the square roots of the eigenvalues of the "interior problem" (in $B(0, R)$, $R \sim R_0$, with Dirichlet boundary conditions) plus the resonances for the "exterior problem", i.e., that in the exterior of the sphere $|x| = R_0$, with Dirichlet boundary conditions (Neumann boundary conditions would not change that). The factor 2 is explained by the fact that each eigenvalue $\lambda^2$ has two square roots: $-\lambda$ and $\lambda$, and resonances are symmetric about $i\mathbf{R}$. In fact, we have a stronger estimate with the extra term $-2(R^n - R_0^n)\tau_n r^n$ in the r.h.s. of the second inequality in the theorem that makes the principal terms independent of $R$. Since the main application of this estimate however is to take th elimit $R \to R_0+$, this term does not add anything new unless one needs sharper estimates on the remainder term.

One can assume more generally that $N^\sharp(r) = O(r^{n^\sharp})$ with some $n^\sharp > n$. Then the result is still true with a remainder $o(r^{n^\sharp})$. Then we recover the asymptotic formula for $N(r)$ in case that $N^\sharp(r)/r^{n^\sharp} \geq 1/C$, see [V5] and [Sj1]: $N(r) = 2N^\sharp(r)(1 + o(r))$. It is known that outside any sector near the real line, the resonances are $O(r^n)$, so the asymptotic is valid actually in any fixed sector around $\mathbf{R}$, as in those works. Also, a "bottle type" theorem can be deduced from (7) but it would not add anything new to the results in [Sj1].

Consider the following examples:

(i) $P_{\mathcal{O}} = -\Delta$ in the domain $\Omega = \mathbf{R}^n \setminus \mathcal{O}$ with Dirichlet or Neumann boundary conditions, where $\mathcal{O} \subset \overline{B(0, R_0)}$, $\partial \mathcal{O} \in C^\infty$;

(ii) $P_V = -\Delta + V(x)$, $V \in L^\infty$, $\text{supp } V \subset \overline{B(0, R_0)}$;

(iii) $P_{g,c} = -c^2 \Delta_g$, where $\Delta_g$ is the Laplace-Beltrami operator associated with a smooth metric $g$, $0 < c = c(x)$ is smooth, such that $g = \delta_{ij}, c = 1$ outside $B(0, R_0)$.

**Corollary 1** *In cases (i), (ii), (iii), one has*

$$M(r) \leq 2A^\sharp r^n + A_{S^{n-1}} R_0^n r^n + o(r^n),$$

*where*

$$A^\sharp = \begin{cases} (2\pi)^{-n} \int_{x \in B(0, R_0) \setminus \mathcal{O}, |\xi| \leq 1} dx\, d\xi, & \text{if } P = P_{\mathcal{O}}, \\ (2\pi)^{-n} \int_{|x| \leq R_0, |\xi| \leq 1} dx\, d\xi, & \text{if } P = P_V, \\ (2\pi)^{-n} \int_{|x| \leq R_0, \sum c^2(x) g^{ij}(x)\xi_i\xi_j \leq 1} dx\, d\xi, & \text{if } P = P_{g,c}. \end{cases}$$

Theorem 1 above shows that the estimate above is sharp in case (i). Next theorem, proven in [Z2] (see the remarks above), shows that our estimate is sharp in case (ii) as well.

**Theorem 3** *Let $V(x) = v(|x|)$ be a radially symmetric potential in $\mathbf{R}^n$, n odd,*

$$v \in C^2([0, R_0]), \quad v(R_0) \neq 0,$$

*and let V be extended as 0 for $|x| > R_0$. Then for the counting function $N(r)$ of $P = -\Delta + V$ we have*

$$N(r) = \left(2\frac{\text{vol}^2(B(0, 1))}{(2\pi)^n} + A_{S^{n-1}}\right) R_0^n r^n + o(r^n), \quad r \to \infty.$$

Finally, we show that the estimate in Corollary 1 is sharp in the "transparent obstacle" case that can be considered as (iii) with singular $c$. Fix $0 < c \neq 1$, and let

$$P = -\tilde{c}^2(x)\Delta, \quad \text{where } \tilde{c}(x) = c \text{ for } |x| \leq R_0, \tilde{c}(x) = 1 \text{ otherwise,}$$

with domain $H^2(\mathbf{R}^n)$ (that corresponds to transmission conditions requiring that $u$ and $\partial u/\partial \nu$ agree on $|x| = R_0$). The operator $P$ is self-adjoint on $L^2(\mathbf{R}^n, \tilde{c}^{-2}dx)$ and satisfies the black-box assumptions. Resonances of $P$ in strips near the real axis for general strictly convex domains and $P$'s of variable coefficients have been studied by Cardoso, Popov and Vodev, see [CPV] and the references there. If $c < 1$, then there are resonances converging rapidly to the real axis; if $c > 1$ there is a resonance free zone $-\Im\lambda \leq (C|\Re\lambda|)^{-1}$, $|\Re\lambda| > C$. This can be explained by the existence, in the case $c < 1$, of totally reflected rays in the interior, close enough to tangent ones to the boundary. In both cases, there is a Weyl type of asymptotic in the strip $0 \leq -\Im\lambda \leq C$ with a suitable $C$. We refer to [CPV] for more results and details. We are concerned here with all the resonances however and we show that for all admissible values of $c$, the estimate in Theorem 2 turns into asymptotic, as $R \to R_0+$.

**Theorem 4** *Let $P$ be the "transparent obstacle" operator as above with some $R_0 > 0$, $c > 0$, $c \neq 1$. Then*

$$N(r) = 2\frac{1}{(2\pi)^n}\int_{|x|\leq R_0; c^2|\xi|^2\leq 1} dx\, d\xi\, r^n + C_{S^{n-1}} R_0^n r^n + o(r^n).$$

## 2 Short review of scattering theory in the black box setting

We introduce briefly the black-box scattering formalism, for more details, see [SjZ1] or [Sj2] for more recent treatment. Fix $R_0 > 0$ and let $\mathcal{H}$ be the complex Hilbert space

$$\mathcal{H} = \mathcal{H}_{R_0} \oplus L^2(\mathbf{R}^n \setminus B(0, R_0)),$$

where $B(0, R_0)$ is the open ball with radius $R_0$ centered at 0. Let $P \geq -C$ be a selfadjoint operator in $\mathcal{H}$ with domain $\mathcal{D}$. Denote by $\mathbf{1}_{B(0,R_0)}$, $\mathbf{1}_{\mathbf{R}^n\setminus B(0,R_0)}$ the corresponding orthogonal projections, and for any $\chi \in L^\infty$ that is constant on $B(0, R_0)$, we define $\chi u$ in an obvious way. In particular, if $\chi_K$ is the characteristic function of $K \supset B(0, R_0)$, we use the notation $\chi_K = \mathbf{1}_K$. Assume that the restriction of $\mathcal{D}$ to $\mathbf{R}^n \setminus B(0, R_0)$ is included in $H^2(\mathbf{R}^n \setminus B(0, R_0))$, and conversely, every $u \in H^2(\mathbf{R}^n \setminus B(0, R_0))$ vanishing near $B(0, R_0)$ belongs to $\mathcal{D}$. The operator $P$ is a compactly supported perturbation of the Laplacian, i.e.,

$$Pu|_{\mathbf{R}^n\setminus B(0,R_0)} = -\Delta u|_{\mathbf{R}^n\setminus B(0,R_0)}.$$

We also require that

$$\mathbf{1}_{B(0,R_0)}(P+i)^{-m_0}$$

to be trace class for some $m_0 > 0$ (see [C]).

We define a reference operator $P^\sharp$ as follows. Fix $R > R_0$ and let $\mathbf{T}_R$ be the flat torus obtained by identifying the opposite sides of $\{x \in \mathbf{R}^n; |x_i| < R, i = 1, \ldots, n\}$. Let $P_{\mathbf{T}}^\sharp$ be the selfadjoint operator defined by

$$P_{\mathbf{T}}^\sharp u = P\chi u + \Delta_{\mathbf{T}_R}(1 - \chi)u, \tag{8}$$

where $\chi = 1$ near $B(0, R_0)$, $\operatorname{supp} \chi \subset B(0, R)$, and $\Delta_{\mathbf{T}_R}$ is the Laplacian on $\mathbf{T}_R$. Then $P_{\mathbf{T}}^\sharp$ is independent on the choice of $\chi$. Our assumptions guarantee that $P_{\mathbf{T}}^\sharp$ has discrete spectrum only, and we set

$$N_{\mathbf{T}}^\sharp(r) = \#\left\{\lambda_j;\ \lambda_j^2 \text{ is an eigenvalue of } P_{\mathbf{T}}^\sharp,\ 0 \leq \lambda_j \leq r\right\}, \tag{9}$$

including multiplicities. Note that $P_{\mathbf{T}}^\sharp$ may have a finite number of negative eigenvalues but they are not included in the counting function above. We assume that

$$N_{\mathbf{T}}^\sharp(r) = O(r^n), \quad N_{\mathbf{T}}^\sharp(r+1) - N_{\mathbf{T}}^\sharp(r) = o(r^n), \quad \text{as } r \to \infty. \tag{10}$$

In most interesting situations, $n^\sharp = n$, and the term $o(r^n)$ can be replaced by $O(r^{n-1})$.

Under the conditions above, T. Christiansen [C] proved that the scattering phase $\sigma(\lambda)$, see (22) below, admits the asymptotic

$$\sigma(r) = N_{\mathbf{T}}^\sharp(r) - \tilde{\tau}_n r^n + o(r^n), \quad \text{as } r \to \infty, \tag{11}$$

where

$$\tilde{\tau}_n = (2\pi)^{-n} \operatorname{vol} B(0,1) \operatorname{vol} \mathbf{T}_R$$

is the Weyl constant related to the torus $\mathbf{T}_R$. As shown in [C], up to $o(r^n)$, $N_{\mathbf{T}}^\sharp(r) - \tau_n r^n$ is independent of the choice of $R > 0$, and in most interesting cases can be expressed by Weyl terms related to $P$ only, see Corollary 1 and its proof. The asymptotics (11) generalizes earlier results in the classical situations and uses techniques developed by Robert [R].

Instead of the reference operator defined above, we consider a reference operator defined in $\mathcal{H}_{R_0} \oplus L^2(B(0,R) \setminus B(0,R_0))$, where $R > R_0$ is fixed. We define $P_{\mathbf{B}}$ to be equal to $P$ on that space and satisfy Dirichlet boundary conditions on $|x| = R$ (in other words, we use an obvious modification of (8)). The results in [C], see Proposition 2.1 there that also holds for manifolds with boundary, imply that $N_{\mathbf{T}}^\sharp(r) - \tilde{\tau}_n r^n = N_{\mathbf{B}}^\sharp(r) - R^n \tau_n + o(r^n)$. From now on, we use $P_{\mathbf{B}}^\sharp$ as a reference operator and will drop the subscript $\mathbf{B}$, i.e., we will denote $P^\sharp = P_{\mathbf{B}}^\sharp$, and $N^\sharp(r)$ is as in (9) but related to $P^\sharp$. Then we have, as in (11),

$$\sigma(r) = N^\sharp(r) - \tau_n r^n + o(r^n), \quad \text{as } r \to \infty. \tag{12}$$

Under the conditions above, $P$ may have a finite number of negative eigenvalues $-\mu_j^2$, and positive eigenvalues as well (the positive ones do not exists in the interesting cases). The resolvent $R(\lambda) = (P - \lambda^2)^{-1} : \mathcal{H}_{\text{comp}} \to \mathcal{H}_{\text{loc}}$ admits a meromorphic continuation from the upper half-plane $\Im \lambda > 0$, where it has poles at $i\mu_j$ only, into the whole complex plane (for $n$ odd), see e.g., [SjZ1] or [Sj2]. We will denote this continuation by $R(\lambda)$. The poles in $\Im \lambda < 0$ are called resonances.

We recall some facts about scattering theory for black boxes, see e.g., [PZ1], [S3] where this is done in the semi-classical setting and we will translate this into the classical setting.

Fix $R_{1,2,3}$ such that $R_0 < R_1 < R_2 < R_3$, and choose a smooth cut-off function $\chi_1$ such that $\chi_1 = 1$ on $B(0, R_1)$, and $\chi_1 = 0$ outside $B(0, R_2)$. For any $\theta \in S^{n-1}$, and any $\lambda > 0$, we are looking for a solution $\psi(x, \theta, \lambda)$ to the problem $(P - \lambda^2)\psi = 0$, $\psi \in \mathcal{D}_{\text{loc}}(P)$ such that

$$\psi = (1 - \chi_1)e^{i\lambda\theta\cdot x} + \psi_{\text{sc}}, \tag{13}$$

with $\psi_{\text{sc}}$ satisfying the Sommerfeld outgoing condition at infinity: $(\partial/\partial r - i\lambda)\psi_{\text{sc}} = O(r^{-(n+1)/2})$, as $r = |x| \to \infty$. Then

$$\psi(x, \theta, \lambda) = e^{i\lambda\theta\cdot x} + \frac{e^{i\lambda r}}{r^{(n-1)/2}} A\left(\frac{x}{r}, \theta, \lambda\right) + O\left(\frac{1}{r^{(n+1)/2}}\right), \quad \text{as } r = |x| \to \infty. \tag{14}$$

The function $A(\omega, \theta, \lambda)$ is the scattering amplitude related to $P$. In order to justify this definition, we will show that $\psi_{\text{sc}}$ is well defined and the limit above exists.

Before proceeding, we will recall the definition for outgoing solution in the case that $\lambda$ is not necessarily real that we will need later. In short, "outgoing" function is a function equal for large $x$ to $R_0(\lambda)f$ for some compactly supported $f$. Here $R_0(\lambda) : \mathcal{H}_{\text{comp}} \to \mathcal{H}_{\text{loc}}$ is the outgoing free resolvent, i.e., the analytic continuation of $R_0(\lambda) = (-\Delta - \lambda^2)^{-1}$ from the upper half-plane into the lower half-plane in $\mathbf{C}$. The extension from the lower to the upper half plane is called incoming.

**Definition 1** *Given $\lambda \in \mathbf{C}$, we say that the function $u$ is $\lambda$-outgoing (or simply, outgoing, if $\lambda$ is understood from the context), if there exists $a > 0$ and $f \in \mathcal{H}_{\text{comp}}$ such that $u|_{|x|>a} = R_0(\lambda)f|_{|x|>a}$.*

Similarly one defines incoming functions.

**Proposition 1 ([S1], see also Lemma 1 in [Z4])**

*(a) For any $f \in \mathcal{H}_{\text{comp}}$ and any $\lambda$ not a resonance, the function $u = R(\lambda)f$ is $\lambda$-outgoing. Moreover, if $\chi$ is a smooth cut-off function such that $\chi = 1$ for $|x| > a$, and $\chi = 0$ in a neighborhood of $B(0, R_0)$ and $\mathrm{supp}\, f$, then we have $R(\lambda)f|_{|x|>a} = -R_0(\lambda)[\Delta, \chi]R(\lambda)f|_{|x|>a}$.*

*(b) Assume $u \in \mathcal{D}_{\text{loc}}(P)$, $(P - \lambda^2)u = f \in \mathcal{H}_{comp}$, $\lambda$ is not a resonance, and $u$ is $\lambda$-outgoing. Then $u = R(\lambda)f$.*

The scattering solution $\psi_{\text{sc}}$ can be constructed as follows. Apply $P - \lambda^2$ to $\psi_{\text{sc}}$ to get

$$(P - \lambda^2)\psi_{\text{sc}} = -(P - \lambda^2)(1 - \chi_1)e^{i\lambda\theta\cdot x} = -[\Delta, \chi_1]e^{i\lambda\theta\cdot x}. \tag{15}$$

Then, since $\psi_{\text{sc}}$ is outgoing, by Proposition 1(b),

$$\psi_{\text{sc}}(x, \theta, \lambda) = -R(\lambda)[\Delta, \chi_1]e^{i\lambda\theta\cdot x}. \tag{16}$$

Choose a smooth function $\chi_2$ with $\mathrm{supp}\,\chi_2 \subset B(0, R_3)$ and $\chi_2 = 1$ on $B(0, R_2) \supset \mathrm{supp}\,\chi_1$. Then, by Proposition 1(a),

$$(1 - \chi_2)\psi_{\text{sc}}(x, \theta, \lambda) = R_0(\lambda)[\Delta, \chi_2]R(\lambda)[\Delta, \chi_1]e^{i\lambda\theta\cdot x},$$

To take the asymptotic as $x = r\omega$, $r = |x| \to \infty$, we recall the asymptotic formula for $R_0(\lambda)f$, where $f$ has compact support, see [M3, section 1.7], (note that in [M3], we have to take complex conjugate since the resonances there are in the upper half-plane)

$$[R_0(\lambda)f](r\omega) = \frac{e^{i\lambda r}}{r^{\frac{n-1}{2}}}\left(v_\infty(\omega, \lambda) + O\left(\frac{1}{r}\right)\right), \tag{17}$$

where

$$v_\infty = \frac{i}{2}(2\pi)^{-\frac{n+1}{2}}\lambda^{\frac{n-3}{2}}e^{-i\pi\frac{n-1}{4}}\hat{f}(\lambda\omega). \tag{18}$$

The function $v_\infty$ is called in the applied literature the far-field pattern of the outgoing solution $v$ to $(-\Delta - \lambda^2)v = 0$ for large $x$ (which always can be expressed as $v = R_0(\lambda)f$ for large $x$). In our case, $v_\infty$ is just the scattering amplitude, if $v = \psi_{\text{sc}}$. Thus we get

$$A(\omega, \theta, \lambda) = \frac{1}{2}e^{-i\pi\frac{n-3}{4}}(2\pi)^{-\frac{n+1}{2}}\lambda^{\frac{n-3}{2}}\int e^{-i\lambda\omega\cdot x}[\Delta, \chi_2]R(\lambda)[\Delta, \chi_1]e^{i\lambda\theta\cdot\bullet}\,dx. \tag{19}$$

It is clear from this formula, that the scattering amplitude $A$ can be extended meromorphically everywhere, where the resolvent admits continuation as well. In particular, all poles of $A$ are poles of the cut-off resolvent are as well.

As in [Z4], [PZ2], introduce the operators

$$[\mathbf{E}_\pm(\lambda)f](\omega) = \int e^{\pm i\lambda\omega\cdot x}f(x)\,dx = \hat{f}(\mp\lambda\omega), \quad \omega \in S^{n-1},$$

and we will apply $\mathbf{E}_\pm(\lambda)$ only to functions $f$ with compact support. Let ${}^t\mathbf{E}_\pm(\lambda)$ be the transpose operators defined as operator with Schwartz kernels ${}^tE(x, \omega) = E(\omega, x)$. Then viewing the scattering amplitude as an operator $\mathbf{A}(\lambda)$ on $L^2(S^{n-1})$ with kernel $A(\omega, \theta, \lambda)$, we recover the formula for $A$ in [PZ2] modulo normalizing factors:

$$\mathbf{A}(\lambda) = \frac{1}{2}e^{-i\pi\frac{n-3}{4}}(2\pi)^{-\frac{n+1}{2}}\lambda^{\frac{n-3}{2}}\mathbf{E}_-(\lambda)[\Delta, \chi_2]R(\lambda)[\Delta, \chi_1]{}^t\mathbf{E}_+(\lambda). \tag{20}$$

The scattering matrix $S(\lambda)$ is an operator on $L^2(S^{n-1})$ and the kernel of $S - I$ is given by

$$a(\omega, \theta, \lambda) := -2\left(\frac{i\lambda}{4\pi}\right)^{\frac{n-1}{2}}A(\omega, \theta, \lambda).$$

Therefore,

$$S(\lambda) = I + c_n\lambda^{n-2}\mathbf{E}_-(\lambda)[\Delta, \chi_2]R(\lambda)[\Delta, \chi_1]{}^t\mathbf{E}_+(\lambda), \quad c_n = -i(2\pi)^{-n}2^{(1-n)/2}. \tag{21}$$

Note that one can replace $\mathbf{E}_-(\lambda)$ by $\mathbf{E}_-(\lambda)\mathbf{1}_{\{R_2<|x|<R_3\}}$, and ${}^t\mathbf{E}_+(\lambda)$ by $\mathbf{1}_{\{R_1<|x|<R_2\}}{}^t\mathbf{E}_+(\lambda)$ above.

The scattering poles (resonances) are defined as the poles of $S(\lambda)$ in the lower halfplane $\Im\lambda < 0$. It is known that $S^*(\lambda) = S^{-1}(\bar\lambda)$, and in particular, $S$ is unitary for $\lambda \in \mathbf{R}$. Note that possible non-negative eigenvalues of $P$ do not contribute to real scattering poles because $\|S(\lambda)\| = 1$ for $\lambda$ not a scattering pole, and if $\lambda_0 \in \mathbf{R}$ were a scattering pole, then we would have $\|S(\lambda)\| \to \infty$, as $\mathbf{R} \ni \lambda \to \lambda_0$. On the other hand, the finite number of negative eigenvalues $-\mu_j^2$ contribute to poles of $S(\lambda)$ at $i\mu_j$ in the physical halfplane $\Im\lambda > 0$ that we do not include in the definition of resonances. It is known, that this definition of resonances is equivalent to the one as the poles of the resolvent in $\Im\lambda < 0$ given above, including the multiplicities.

The *scattering determinant* $s(\lambda)$ is defined by

$$s(\lambda) = \det S(\lambda).$$

The *scattering phase* $\sigma(\lambda)$ is given by

$$\sigma(\lambda) = \frac{1}{2\pi\mathrm{i}} \log s(\lambda), \quad \sigma(0) = 0, \quad \sigma(-\lambda) = -\sigma(\lambda). \tag{22}$$

## 3 Preliminary results

**Lemma 1** *Let $M(r)$ be as in (4). Then*

$$M(r) = Ar^n + o(r^n), \quad \text{as } r \to \infty. \tag{23}$$

*with some $A > 0$, if and only if*

$$N(r) = Ar^n + o(r^n), \quad \text{as } r \to \infty, \tag{24}$$

*Proof:* Assume (23), i.e., $M(r) = Ar^n + \mu(r)r^n$, where $\lim_{r\to\infty} \mu(r) = 0$. Set $\mu_+(r) = \sup_{t \geq r} |\mu(t)|$. Then $|M(r) - Ar^n| \leq \mu_+(r)r^n$, and $\mu_+$ is decreasing and converges to 0. If $\mu_+(r) = 0$ for $r$ large enough, our statement follows easily. Assume this never happens, then $\mu_+(r) > 0$ for all $r$. Set $\alpha = r\sqrt{\mu_+(r)}$. Then

$$n\int_r^{r+\alpha} \frac{N(t)}{t}\,\mathrm{d}t = M(r+\alpha) - M(r) = An\alpha r^{n-1} + O(r^n \mu_+(r)).$$

On the other hand,

$$n\alpha \frac{N(r)}{r+\alpha} \leq n\int_r^{r+\alpha} \frac{N(t)}{t}\,\mathrm{d}t \leq n\alpha \frac{N(r+\alpha)}{r}.$$

Therefore,

$$N(r) \leq Ar^{n-1}(r+\alpha) + C'r^n\mu_+(r)\left(\frac{r}{\alpha(r)} + 1\right) = Ar^n + o(r^n).$$

Similarly,

$$N(r+\alpha) \geq Ar^n - C''r^{n+1}\mu_+(r)/\alpha = Ar^n - C''r^n\sqrt{\mu_+(r)}.$$

Replace $r$ by $r - \alpha(r)$ to finish the proof of the implication (23) $\implies$ (24).

Assume now (24). Given $\epsilon > 0$, let $a$ be such that $|N(r) - Ar^n| \leq \epsilon r^n$ for $r > a$. Then

$$|M(r) - Ar^n| \leq n\int_0^r \frac{|N(t) - At^n|}{t}\,\mathrm{d}t \leq C(a) + \epsilon\int_a^r nt^{n-1}\,\mathrm{d}t = C(a) + \epsilon(r^n - a^n).$$

Divide by $r^n$ to get $r^{-n}|M(r) - Ar^n| \leq 2\epsilon$ for $r$ large enough, and this proves (23). $\square$

The following lemma is due essentially by R. Froese [Fr] and its semiclassical version is presented in [PZ2].

**Lemma 2** *For any $r > 0$ we have*

$$\frac{1}{n}M(r) = 2\int_0^r \frac{\sigma(t)}{t}\,dt + \frac{1}{2\pi}\int_0^\pi \log|s(re^{i\theta})|\,d\theta + m(r),$$

*where $0 \leq m(r) = O(\log r)$ (and $m = 0$ if P has no negative eigenvalues).*

*Proof:* The resonances are zeros of $s(\lambda)$ in $\Im\lambda > 0$, with multiplicities with finitely many possible exceptions at points in the set $\{i\mu_j\}$. On the other hand, $s(\lambda)$ may have a finite number of poles in the same set. Assume first that $r$ is not an absolute value of a resonance or a zero of $s(\lambda)$ in $\Im\lambda > 0$. Let $n(t)$ be the number of poles of $s(\lambda)$ on $i(0,t)$. Following the proof of Jensen's formula, we integrate $s'/s$ along the contour $[-r, r] \cup r\exp(i[0,\pi])$ keeping in mind that $s'/s = 2\pi i\sigma'$, to get

$$\begin{aligned} N(t) - n(t) &= \frac{1}{2\pi i}\oint \frac{s'(z)}{s(z)}\,dz = \Im\frac{1}{2\pi}\oint \frac{s'(z)}{s(z)}\,dz \\ &= \int_{-t}^t \sigma'(z)\,dz + \frac{1}{2\pi}\int_0^\pi t\frac{d}{dt}\log|s(te^{i\theta})|\,d\theta \\ &= 2\sigma(t) + \frac{1}{2\pi}\int_0^\pi t\frac{d}{dt}\log|s(te^{i\theta})|\,d\theta. \end{aligned}$$

Divide by $t$ and integrate to get the lemma. Note that the integrand has singularities at the resonances and the zeros, and to justify the calculations we use the same arguments as in [T] together with the fact that $s(0) = 1$. □

**Proposition 2**

$$M(r) = 2\left(N^\sharp(r) - \tau_n R^n r^n\right) + \frac{n}{2\pi}\int_0^\pi \log|s(re^{i\theta})|\,d\theta + o(r^n), \quad \text{as } r \to \infty. \tag{25}$$

*Proof:* We apply Lemma 2. To estimate the scattering phase, we apply the asymptotic (12). □

Below, we will sketch a proof that the integral term in (25) is bounded by $C(R_0^n r^n + 1)$ with an absolute constant $C$, which is one of the ways to prove the polynomial bound (1) in the general case if $n^\sharp = n$, with obvious modifications if $n^\sharp > n$, see e.g., [PZ2]. The reason we sketch this proof is to explain the main idea in the proof of Theorem 2.

To estimate the scattering determinant $s(\lambda)$, we proceed in the usual way, see for example [PZ2]. By (21), we need to estimate the characteristic values of $S(\lambda) - I$, which equals $\mathbf{A}(\lambda)$ modulo polynomial factors, see (20) and (21). This reduces to an estimate of the characteristic values of the operators $\mathbf{E}_-(\lambda)\mathbf{1}_{R_2 < |x| < R_3}$, and $\mathbf{1}_{R_1 < |x| < R_2}\mathbf{E}_+(\lambda)$, see the remark after (21), and the latter can be done by estimating

$$\Delta_\theta^m e^{i\lambda x\cdot\theta}, \quad |x| \leq R_3.$$

We need to work here in a sector $0 < \delta \leq \arg\lambda \leq \pi - \delta$. Using a standard argument, to cover the missing sectors $0 \leq \arg\lambda \leq \delta$ and $\pi - \delta \leq \arg\lambda \leq \pi$, we use the fact that $|s| = 1$ on $\mathbf{R}$ and the Phragmén-Lindelöf principle.

More precisely, assume that $R_3 = 1$ and that in the representation (21), the cut-off functions $\chi_1$ and $\chi_2$ are so that they are supported in $B(0,1)$. Then the statement for any $R_3 > 0$ would follow by a scaling argument. In what follows, $0 < \delta \leq \arg\lambda \leq \pi - \delta$ for some $0 < \delta < 1/(n+1)$ and $|\lambda| \geq 1$. Note that $\|R(\lambda)\| \leq C/|\lambda|^2 \leq C$ for $\lambda$ in this region, so in (21), we have

$$\left\|[\Delta, \chi_2]R(\lambda)[\Delta, \chi_1]e^{i\lambda\theta\cdot\bullet}\right\|_{L^2(\mathbf{R}^n)} \leq Ce^{|\lambda|} \tag{26}$$

with $C > 0$ depending on $\delta$ only and in particular, independent of $P$. Therefore, for any $m = 1, 2, \ldots$,

$$|\Delta_\omega^m a(\omega, \theta, \lambda)| \leq Ce^{|\lambda|} \max_{|x|\leq 1}\left|\Delta_\omega^m e^{-i\lambda x\cdot\omega}\right|, \quad \forall\omega, \theta,$$

with a similar $C$. This shows that one can get the standard now estimate of the characteristic values of the operator with kernel $a(\omega, \theta, \lambda)$ as in [Z2, Lemma 2]

$$\left|\Delta_\omega^m a(\omega, \theta, \lambda)\right| \leq C^{2m+1} \left(|\lambda|^{2m} + (2m)!\right) e^{2|\lambda|},$$

with $C$ as above. By [Z2, Proposition 2],

$$|s(\lambda)| \leq C e^{C|\lambda|^n} \tag{27}$$

with $C > 0$ independent of $P$ under the assumption $R_0 < R_3 = 1$. This is an analogue of [PZ2, Lemma 4.3], where (27) is proved in the semiclassical case (and is implicit in [PZ1]). A scaling argument gives us immediately

$$|s(\lambda)| \leq C e^{C R_3^n |\lambda|^n}, \quad C = C(n, R_3/R_0),$$

where $R_3$ is any constant such that $R_3 > R_0$. As mentioned above, using Phragmén-Lindelöf principle, we extend this to $\Im \lambda \geq 0$.

The main idea behind the proof of Theorem 2 is the following. To get an explicit value for $C$, we notice that by a well known expansion of $e^{i\lambda x \cdot \theta}$ in spherical harmonics and Bessel functions $J_\nu$, see Lemma 3, one can find the characteristic values of $\mathbf{E}_-(\lambda) \mathbf{1}_{R_2 < |x| < R_3}$, and $\mathbf{1}_{R_1 < |x| < R_2} \mathbf{E}_+(\lambda)$ explicitly in terms of $J_\nu(r)$, see (46). In the spherically symmetric cases (i), (ii), (iii), this in fact gives not only an upper bound, but an asymptotic of the integral term in (25).

# 4  Preliminaries about Bessel's functions

We will recall some facts about separation of variables in polar coordinates for the Laplace operator, see e.g., [Fo], and some asymptotics of Bessel's functions, see [O1], [O2], [O3]. Denote by $Y_l^m$, $l = 0, 1, \ldots$, $m = 1, \ldots m(l)$, an orthonormal set of spherical harmonics on $S^{n-1}$. They are the eigenfunctions of the Laplacian $\Delta_{S^{n-1}}$ on $S^{n-1}$. We have

$$-\Delta_{S^{n-1}} Y_l^m = l(l + n - 2) Y_l^m, \quad l = 0, 1, \ldots; \ m = 1, \ldots, m(l).$$

For each $l$, the multiplicity of the eigenvalue $l(l+n-2)$ is given by

$$m(l) = \frac{2l + n - 2}{n - 2} \binom{l + n - 3}{n - 3} = \frac{2l^{n-2}}{(n-2)!} \left(1 + O(l^{-1})\right). \tag{28}$$

Any solution $u$ of the Helmholtz equation $(-\Delta - \lambda^2) u = 0$ near 0 has the form

$$u(x) = \sum_{l=0}^\infty c_{lm} (\lambda r)^{1-n/2} J_{l+n/2-1}(\lambda r) Y_l^m(\omega), \tag{29}$$

where $x = r\omega$ and $r > 0$, $|\omega| = 1$ are polar coordinates. Similarly, any outgoing solution at $\infty$ has similar expansion, with $J_\nu$ replaced by $H_\nu^{(2)}$. The functions $\lambda^{1-n/2} J_{l+n/2-1}(\lambda)$ are entire and in particular, regular at $\lambda = 0$.

We will need the formula below.

**Lemma 3** *For any $\theta \in S^{n-1}$, $\lambda \in \mathbf{C}$, and $x \in \mathbf{R}^n$, we have*

$$e^{i\lambda x \cdot \theta} = (2\pi)^{n/2} \sum_l i^l Y_l^m(\omega) \overline{Y_l^m}(\theta) (\lambda r)^{1-n/2} J_{l+n/2-1}(\lambda r), \quad x = r\omega. \tag{30}$$

*Proof:* This formula is known and widely used, at least in the 3D case. We could not find a proof for general odd $n$'s, so we will sketch one here.

Note first that the series above converges absolutely and uniformly for any $\theta$, and $\lambda$ in any compact, as a consequence of the well known asymptotics of $J_\nu$, as $\nu \to \infty$. It is enough to prove it for real $\lambda$, because we can then extend it analytically for all $\lambda$. With $x = r\omega$, we have

$$e^{i\lambda r \omega \cdot \theta} = \sum_{l,m} a_{lm} Y_l^m(\omega),$$

where
$$a_{lm} = \int_{S^{n-1}} e^{i\lambda r \omega \cdot \theta} \overline{Y_l^m}(\omega) \, d\omega.$$

By the Funk-Hecke Theorem (see e.g., [EMOT, §11.4]),

$$a_{lm} = \overline{Y_l^m}(\theta) \, i^l \, (2\pi)^{n/2} \int t^{1-n/2} J_{l+n/2-1}(t) f(t) \, dt, \tag{31}$$

where

$$f(t) = (2\pi)^{-1} \int_{-1}^{1} e^{-ist} e^{i\lambda rs} \, ds. \tag{32}$$

The well-known integral representation

$$\Gamma\left(\nu + \frac{1}{2}\right) J_\nu(z) = \frac{1}{\sqrt{\pi}} \left(\frac{z}{2}\right)^\nu \int_{-1}^{1} e^{izt} \left(1 - t^2\right)^{\nu - 1/2} dt$$

shows that $t^{1-n/2} J_{l+n/2-1}(t)$ has Fourier transform supported in $[-1, 1]$. Expressing (31) via the Plancherel theorem, we see that one can change the definition of $f$ in (32) by integrating from $-\infty$ to $\infty$, and this would not change (31). The integral in (32) however, over the whole real line, is simply $\delta(t - \lambda r)$. Setting $f = \delta(t - \lambda r)$ in (31), we complete the proof of the lemma. $\square$

Let Ai be the Airy function, having its zeros on the negative real axis; set $\text{Ai}_\pm(w) = \text{Ai}(e^{\mp 2\pi i/3} w)$. Then

$$\text{Ai}(w) \sim \frac{e^{-\xi}}{2\sqrt{\pi} w^{1/4}} \left(1 + \sum_{s=1}^{\infty} \frac{c_s}{\xi^s}\right), \quad \text{Ai}'(w) \sim \frac{-w^{1/4} e^{-\xi}}{2\sqrt{\pi}} \left(1 + \sum_{s=1}^{\infty} \frac{c_s'}{\xi^s}\right), \quad |\arg w| < \pi, \tag{33}$$

where $\xi = \frac{2}{3} w^{3/2}$. In particular, $\text{Ai}(\omega)$ is exponentially decaying, as $|\omega| \to \infty$, in the sector $|\arg w| < \pi/3$. The expansions above hold for $\text{Ai}_\pm$ as well with the appropriate choice of the branch of $\omega^{3/2}$; this branch is uniquely determined by the condition that $\text{Ai}_\pm(\omega)$ is exponentially decaying for $\pm \arg w \in (\pi/3, \pi)$.

Near the zeros of $\text{Ai}(w)$ we have [O3, p. 413],

$$\text{Ai}(-w) \sim \frac{1}{\sqrt{\pi} w^{1/4}} \left\{ \cos\left(\xi - \frac{\pi}{4}\right) \left(1 + \sum_{s=1}^{\infty} \frac{d_s}{\xi^{2s}}\right) + \sin\left(\xi - \frac{\pi}{4}\right) \sum_{s=0}^{\infty} \frac{\tilde{d}_s}{\xi^{2s+1}} \right\}, \quad |\arg w| < \frac{2\pi}{3}, \tag{34}$$

$$\text{Ai}'(-w) \sim \frac{w^{1/4}}{\sqrt{\pi}} \left\{ \sin\left(\xi - \frac{\pi}{4}\right) \left(-1 + \sum_{s=1}^{\infty} \frac{d_s'}{\xi^{2s}}\right) + \cos\left(\xi - \frac{\pi}{4}\right) \sum_{s=0}^{\infty} \frac{\tilde{d}_s'}{\xi^{2s+1}} \right\}, \quad |\arg w| < \frac{2\pi}{3}. \tag{35}$$

Following Olver [O1], [O2], introduce the functions

$$\rho(z) = \frac{2}{3} \zeta^{3/2} = \log \frac{1 + \sqrt{1 - z^2}}{z} - \sqrt{1 - z^2}, \quad |\arg z| < \pi.$$

The branches of the functions appearing above are chosen so that $\zeta$ is real, if $z$ is real. The mapping properties of $\rho$ and $\zeta$ can be found in [O3, p. 336], and they are of fundamental importance in our analysis. An important role is played by the eye-shaped domain **K**, symmetric about the real axis, bounded by the following curve and its conjugate:

$$z = \pm (t \coth t - t^2)^{1/2} + i(t^2 - t \tanh t)^{1/2}, \quad 0 \le t \le t_0, \tag{36}$$

and $t_0 = 1.19967864\ldots$ is the positive root of $t = \coth t$. The intercepts of $\partial \mathbf{K}$ with the imaginary axis are $\pm (t_0^2 - 1)^{1/2} = \pm i \, 0.6627\ldots$. Notice that in $\mathbf{C}_+ = \{z; \Im z > 0\}$, we have $\Re \rho > 0$ in **K**, and $\Re \rho < 0$ outside $\bar{\mathbf{K}}$.

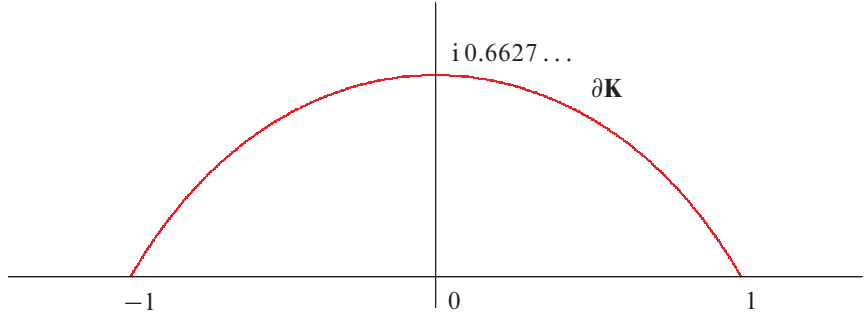

Figure 1: Sketch of the domain **K** in the upper halfplane $\Im z \geq 0$.

The following asymptotic expansions are established in [O1], [O2], see [O1, Theorem B], and [O2, §4]:

$$J_\nu(\nu z) \sim \left(\frac{4\zeta}{1-z^2}\right)^{1/4} \left(\frac{\text{Ai}(\nu^{2/3}\zeta)}{\nu^{1/3}} \sum_{s=0}^\infty \frac{A_s(\zeta)}{\nu^{2s}} + \frac{\text{Ai}'(\nu^{2/3}\zeta)}{\nu^{5/3}} \sum_{s=0}^\infty \frac{B_s(\zeta)}{\nu^{2s}}\right), \tag{37}$$

$$H_\nu^{(1,2)}(\nu z) \sim \frac{2e^{\mp i\pi/3}}{\nu^{1/3}} \left(\frac{4\zeta}{1-z^2}\right)^{1/4} \left(\frac{\text{Ai}_\mp(\nu^{2/3}\zeta)}{\nu^{1/3}} \sum_{s=0}^\infty \frac{A_s(\zeta)}{\nu^{2s}} + \frac{\text{Ai}'_\mp(\nu^{2/3}\zeta)}{\nu^{5/3}} \sum_{s=0}^\infty \frac{B_s(\zeta)}{\nu^{2s}}\right). \tag{38}$$

The infinite series expansions above are uniform in $|\arg z| \leq \pi - \delta$, $\delta > 0$ fixed. Similar expansions hold for the derivatives, and they can be obtained by differentiating (37), (38) term by term.

## 5 Sharp estimate of the scattering determinant $s(\lambda)$

To prove Theorem 2, we will use Proposition 2. To this end, we need to estimate the integral

$$\frac{n}{2\pi} \int_0^\pi \log|s(re^{i\theta})| \, d\theta. \tag{39}$$

We will prove first the following sharp estimate on the scattering determinant $s(\lambda)$ in $\Im \lambda \geq 0$:

**Theorem 5**
(a) For any $\theta \in [0, \pi]$,
$$\log|s(re^{i\theta})| \leq h_n(\theta) R_0^n r^n + o(r^n), \quad \text{as } r \to \infty, \tag{40}$$

where
$$h_n(\theta) = \frac{4}{(n-2)!} \int_0^\infty \frac{[-\Re\rho]_+(te^{i\theta})}{t^{n+1}} \, dt \tag{41}$$

and the remainder term depends on the operator $P$, and is uniform for $0 < \delta \leq \theta \leq \pi - \delta$ or any $\delta \in (0, \pi)$.
(b) For any $\delta > 0$,
$$\log|s(re^{i\theta})| \leq \left(h_n(\theta) R_0^n + \delta\right) r^n + o(r^n), \quad \text{as } r \to \infty \tag{42}$$

uniformly in $\theta \in [0, \pi]$.

**Remark.** The integral above can be evaluated to some extent. In the $n = 3$ case, for example, we get

$$h_3(\theta) = 4 \int \frac{[-\Re\rho]_+(te^{i\theta})}{t^4} \, dt = \frac{4\Re(1-z^2)^{3/2}}{9|z|^3},$$

where $z$ is the unique point on $\partial \mathbf{K}$ with argument $\theta$, i.e., $z$ is given by (36) with $t \in [0, t_0]$ the unique solution of

$$\tan^2 \theta = \frac{t - \tanh t}{\coth t - t}.$$

Another way to define $z = z(\theta)$ is as the solution of $\Re\rho(z) = 0$, $\arg z = \theta$.

**Remark.** One can verify that $h_n(\theta) \sim C_n \theta^{5/2}$ as $\theta \to 0+$. We prove in fact, that one can replace $\delta$ in (42) by $O(\theta)$ near $\theta = 0$, see (55). A more careful analysis of the leading term in (42) as $\theta \to 0+$ is in principle possible, but not needed for our purposes because at the end, we will integrate (42) in $\theta \in [0, \pi]$.

*Proof of Theorem 5:* We will estimate the integral (39). Recall Weyl's estimate

$$|\det(I + B)| \leq \prod_{j=1}^{\infty} (1 + \mu_j(B)),$$

provided that $B$ is trace class. We also recall that $\det(I + AB) = \det(I + BA)$. Then by (21),

$$\begin{aligned}
\log|s(\lambda)| &= \log\left|\det\left(I + c_n \lambda^{n-2} \mathbf{E}_-(\lambda) \mathbf{1}_{R_2 < |x| < R_3} [\Delta, \chi_2] R(\lambda) [\Delta, \chi_1] \mathbf{1}_{R_1 < |x| < R_2} {}^t\mathbf{E}_+(\lambda)\right)\right| \\
&= \log\left|\det\left(I + c_n \lambda^{n-2} [\Delta, \chi_2] R(\lambda) [\Delta, \chi_1] \mathbf{1}_{R_2 < |x| < R_3} {}^t\mathbf{E}_+(\lambda) \mathbf{E}_-(\lambda) \mathbf{1}_{R_1 < |x| < R_2}\right)\right| \\
&\leq \sum_{j=1}^{\infty} \log\left(1 + \mu_j\left(c_n \lambda^{n-2} [\Delta, \chi_2] R(\lambda) [\Delta, \chi_1] \mathbf{1}_{R_2 < |x| < R_3} {}^t\mathbf{E}_+(\lambda) \mathbf{E}_-(\lambda) \mathbf{1}_{R_1 < |x| < R_2}\right)\right). \quad (43)
\end{aligned}$$

We work in the set

$$\epsilon_0 \leq \arg \lambda \leq \pi - \epsilon_0, \quad 2 \leq |\lambda| \quad (44)$$

with a fixed $0 < \epsilon_0 < \pi/2$. There, we have by the spectral theorem and standard elliptic estimates for $-\Delta$,

$$\|[\Delta, \chi_2] R(\lambda) [\Delta, \chi_1]\| \leq C.$$

Here and below, all constants may depend on $\epsilon_0$ that is kept fixed. Use this and $\mu_j(AB) \leq \|A\| \mu_j(B)$ to get

$$\log|s(\lambda)| \leq \sum_{j=1}^{\infty} \log\left(1 + C_0 |\lambda|^{n-2} \mu_j\left(\mathbf{1}_{R_2 < |x| < R_3} {}^t\mathbf{E}_+(\lambda) \mathbf{E}_-(\lambda) \mathbf{1}_{R_1 < |x| < R_2}\right)\right) \quad (45)$$

By (30), the operator

$$L^2(S^{n-1}) \ni f(\omega) \longmapsto \int_{S^{n-1}} e^{\mp i\lambda r \omega \cdot \theta} f(\omega) \, d\omega \in L^2(S^{n-1}), \quad |\theta| = 1,$$

is a diagonal one in the spherical harmonics base, and has eigenvalues $(2\pi)^{n/2} i^l (\lambda r)^{1-n/2} J_{l+n/2-1}(\mp r\lambda)$ with multiplicities $m(l)$ given by (28). Therefore, the non-zero characteristic values $\mu_j\left(\mathbf{1}_{R_2 < |x| < R_3} {}^t\mathbf{E}_+(\lambda) \mathbf{E}_-(\lambda) \mathbf{1}_{R_1 < |x| < R_2}\right)$ coincide with

$$\tilde{\mu}_l = (2\pi)^n \left(\int_{R_1}^{R_2} |(\lambda r)^{1-n/2} J_{l+n/2-1}(-\lambda r)|^2 r^{n-1} \, dr\right)^{1/2} \left(\int_{R_2}^{R_3} |(\lambda r)^{1-n/2} J_{l+n/2-1}(\lambda r)|^2 r^{n-1} \, dr\right)^{1/2}, \quad (46)$$

$l = 0, 1, \ldots$, each one repeated $m(l)$ times. The sequence above may not be decreasing but since the series (45) converges absolutely, it will not be affected by rearrangement of its terms.

So the problem is reduced to that of estimating the exponential growth of

$$|J_\nu(\lambda R_1) J_\nu(-\lambda R_2)|, \quad R_1 \sim R_2 \sim R.$$

Notice first, that by (37),

$$|J_\nu(\nu z)| \leq C e^{-\nu \Re \rho}$$

for $z$ in the sector $\epsilon_0 \leq \arg z \leq \pi - \epsilon_0$. From now on, we denote

$$\nu = l + \frac{n}{2} - 1. \tag{47}$$

Note that $\nu$ is half-integer, because $n$ is odd. Then (see (46)),

$$|(\lambda r)^{1-n/2} J_\nu(\lambda r)|^2 r^{n-1} \leq C|\lambda|^{2-n} r e^{-2\nu \Re \rho(\lambda r/\nu)}.$$

We want to estimate this for $R_1 \leq r \leq R_2$ and $\lambda$ as in (44). Observe that for $t > 0$, $-\mathrm{d}\Re\rho(tz)/\mathrm{d}t = t^{-1} \Re \sqrt{1-(tz)^2} > 0$ for $z$ in the sector (44). Therefore, the exponent above is an increasing function of $r$, and

$$\log \left( |(\lambda r)^{1-n/2} J_\nu(\lambda r)|^2 r^{n-1} \right) \leq -2\nu \Re \rho(\lambda R_2/\nu) + C, \quad R_1 \leq r \leq R_2.$$

This yields (see (46), (45)),

$$\log \left( C_0 |\lambda|^{n-2} \tilde\mu_l \right) \leq -2\nu \Re \rho(\lambda R_3/\nu) + C \log|\lambda|, \tag{48}$$

and

$$\log |s(\lambda)| \leq \sum_{l=0}^\infty m(l) \log \left( 1 + C_0 |\lambda|^{n-2} \tilde\mu_l \right). \tag{49}$$

Observe that $-\Re\rho(\lambda R_3/\nu) > 0$ for $\lambda R_3/\nu \notin \mathbf{K}$, and $-\Re\rho(\lambda R_3/\nu) \leq 0$ otherwise (in the set (44)). Fix $0 < \delta < 1$. We will split the sum above into three parts $\Sigma_{1,2,3}$ corresponding to $\lambda R_3/\nu \notin (1+\delta)\mathbf{K}$, $\lambda R_3/\nu \in (1+\delta)\mathbf{K} \setminus (1-\delta)\mathbf{K}$, and $\lambda R_3/\nu \in (1-\delta)\mathbf{K}$, respectively.

To estimate $\Sigma_1$, we use the inequality

$$0 < a < A, \quad 1 \leq A \implies \log(1+a) \leq \log(1+A) = \log A + \log(1+1/A) \leq \log A + \log 2, \tag{50}$$

and (48) to get

$$\Sigma_1 \leq \sum_{\lambda R_3/\nu \notin (1+\delta)\mathbf{K}} 2\nu m(l)[-\Re\rho](\lambda R_3/\nu) + C|\lambda|^{n-1} \log|\lambda|$$

$$\leq \sum_{\lambda R_3/\nu \notin (1+\delta)\mathbf{K}} \frac{4\nu^{n-1}}{(n-2)!}[-\Re\rho](\lambda R_3/\nu) + C|\lambda|^{n-1} \log|\lambda|.$$

To be more clear, the summation above is taken over all $\nu = l+n/2-1$, $l = 0,1,\ldots$, with the property indicated. To estimate the remainder we used (28) and the fact that $\nu \leq C|\lambda|$. To get the second inequality above, we used (28) again and the fact that $-\Re\rho(z) \leq |z|(1+O(|z|^{-1}))$, $|z| > 1/C$, $\Im z > 0$. The fact that $-\Re\rho(\lambda R_3/\nu)$ is a decreasing function of $\nu$ makes it easy to prove that one can replace the sum above by an integral with the same remainder estimate:

$$\Sigma_1 \leq \int_{\lambda R_3/\nu \notin (1+\delta)\mathbf{K}} \frac{4\nu^{n-1}}{(n-2)!}[-\Re\rho](\lambda R_3/\nu) \, \mathrm{d}\nu + C|\lambda|^{n-1} \log|\lambda|$$

$$= (rR_3)^n \int_{te^{i\theta} \notin (1+\delta)\mathbf{K}} \frac{4}{(n-2)!} \frac{[-\Re\rho](te^{i\theta})}{t^{n+1}} \, \mathrm{d}t + Cr^{n-1} \log r, \quad \lambda := re^{i\theta}. \tag{51}$$

We made the change of variables $rR_3/\nu = t$ above.

To estimate $\Sigma_2$, note that $-\Re\rho(R_3\lambda/\nu)$ changes sign now but we have $-\Re\rho(R_3\lambda/\nu) \leq C\delta$ with $C$ independent of $\epsilon_0$ as can be seen using the relationship $\rho = (2/3)\zeta^{3/2}$ between $\rho$ and $\zeta$, and the fact that $\zeta$ is analytic near $z = 1$. To apply the same argument near $z = -1$, we use the symmetry property $\Re\rho(-\bar{z}) = \Re\rho(z)$. Reasoning as above, we get

$$\Sigma_2 \leq C\delta r^n, \tag{52}$$

uniformly in $\theta$, where $r = |\lambda|$ as above.

Finally, for $\nu$ as in the definition of $\Sigma_3$, we have $\delta/C \leq \Re\rho(R_3\lambda/\nu)$. Using the inequality $\log(1+a) \leq a$, $a > 0$, and (48), we deduce

$$\Sigma_3 \leq Ce^{-r/C}, \tag{53}$$

with some $C = C(\delta)$.

Recall that in $\Im z > 0$, $\mathbf{K}$ is characterized by the condition $\Re\rho(z) > 0$. We can integrate over $te^{i\theta} \notin \mathbf{K}$ in (51) and this will only increase the integral (by adding $O(\delta r^n)$, as in (52), so we are in no danger of losing the sharpness of the estimate). Also, we can integrate over $\Im te^{i\theta} > 0$, i.e., over $t > 0$, as long as we replace $-\Re\rho$ by $[-\Re\rho]_+$ because $[-\Re\rho]_+ = 0$ in $\mathbf{K}$.

Combine (51), (52) and (53) to get

$$\log|s(re^{i\theta})| \leq (rR_3)^n h_n(\theta) + C\delta r^n + C_1 r^{n-1}\log r, \quad \epsilon_0 \leq \theta \leq \pi - \epsilon_0, \tag{54}$$

where $C > 0$ depends on $\epsilon_0$ and $R_3$; $C_1$ depends in addition on $\delta$, see (53), but they not depend on $r$. Clearly, $h_n(\theta) \geq 1/C$ for $\theta$ as in (54), therefore the $C\delta r^n$ term can be absorbed by the preceeding one at the expense of parturbing $R_3$ by an $O(\delta)$ term. So we may formally assume that $\delta = 0$, and (54) still holds.

Consider the function $f_\theta(r) = r^{-n}\log|s(re^{i\theta})| - R_0^n h_n(\theta)$. By (54) and the remark above, for any $R_3 > R_0$, we have $f_\theta(r) \leq C_0(R_3^n - R_0^n) + \alpha_{\theta,R_3}(r)$, where $0 \leq \alpha_{\theta,R_3}(r) \to 0$, as $r \to \infty$ uniformly for $\theta$ as in (54). Set

$$\beta(r) = \inf_{\epsilon_0 \leq \theta \leq \pi - \epsilon_0, R_0 < R_3 < R_0 + 1} \left\{ C_0(R_3^n - R_0^n) + \alpha_{\theta,R_3}(r) \right\}.$$

Choose $\epsilon > 0$. Let $\tilde{R}_3 \in (R_0, R_0 + 1)$ be such that $C_0(\tilde{R}_3^n - R_0^n) \leq \epsilon/2$. Let $r_0$ be such that $\alpha_{\theta,\tilde{R}_3}(r) \leq \epsilon/2$ for $r \geq r_0$ and all $\theta$ as above. Then $\beta(r) \leq \epsilon$ for $r \geq r_0$. In other words, $\beta(r) \to \infty$, as $r \to \infty$, and $f_\theta(r) \leq \beta(r)$ for all $\theta$ as above. This completes the proof of (a).

To prove (b), notice that since the scattering operator is unitary for real $\lambda$, we have $\log|s(\pm r)| = 0$. Also, we have the rather rough a priori estimate $|s(\lambda)| \leq Ce^{C|\lambda|^n}$ in $\Im\lambda > 0$, see e.g., [PZ2]. Since the function $r^n \sin(n\theta)$ is harmonic in $\lambda = re^{i\theta}$, we can apply the Phragmén-Lindelöf principle to the harmonic function $\log|s(re^{i\theta})| - Ar^n\sin(n\theta)$ in the sector $0 \leq \theta \leq \epsilon_0$ with $A = 2R_0^n h_n(\epsilon_0)/\sin(n\epsilon_0)$ to get

$$\log|s(re^{i\theta})| \leq 2R_0^n h_n(\epsilon_0) \frac{\sin(n\theta)}{\sin(n\epsilon_0)} r^n, \quad \text{for } r \gg 1, 0 \leq \theta \leq \epsilon_0. \tag{55}$$

Since $h_n(\theta) \to 0$, as $\theta \to 0$, this proves part (b) of the Theorem. $\square$

Set

$$\partial\mathbf{K}_+ := \partial\mathbf{K} \cap \{\Im z \geq 0\}. \tag{56}$$

**Lemma 4**

$$A_{S^{n-1}} = \frac{2}{\pi} \frac{1}{n(n-2)!} \int_{\partial\mathbf{K}_+} \frac{|1-z^2|^{1/2}}{|z|^{n+1}} |dz|.$$

*Proof:* We start with formula (5). Set $u = -\Re\rho$. Then $u$ is a $C^1$ function in the closure of $\mathbf{C}_+ \setminus \mathbf{K}$. Indeed, this is true at $z = 1$ because $\rho = 2\zeta^{3/2}/3$, and $\zeta$ is analytic there. It is also true at $(-\infty, -1]$ as well, because $u(-\bar{z}) = u(z)$. Inside that domain, $u$ is harmonic. On the other hand, $\Delta n^{-2}|z|^{-n} = |z|^{-n-2}$, where $\Delta$ is the Laplacian in $\mathbf{R}^2$ that we identify with $\mathbf{C}^2$. Now, we view the integrand in (5) as $u\Delta n^{-2}|z|^{-n}$, and apply Green's formula in $\mathbf{C}_+ \setminus \mathbf{K}$ using the fact that $u = O(|z|)$, $\nabla u = O(1)$, as $|z| \to \infty$. Since $u = 0$ on the boundary, we get

$$\int_{\Im z > 0} \frac{[-\Re\rho]_+(z)}{|z|^{n+2}}\, dx\, dy = \int_{\partial \mathbf{K}_+} \frac{\partial u}{\partial \nu} \frac{1}{n^2|z|^n}\, |dz| + 2\int_1^\infty \frac{\partial u}{\partial y} \frac{1}{n^2 x^n}\, dx,$$

where $\nu$ is the unit normal to $\partial \mathbf{K}$, pointing into the exterior to $\mathbf{K}$. Note that in both integrals in the r.h.s. above, the integration is taken over a curve defined by $u = 0$. Then $\nu = \nabla u/|\nabla u|$, and $\partial u/\partial \nu = |\nabla u|$. On the other hand, $\rho' = -\sqrt{1 - z^2}/z$, therefore $|\nabla u| = |1 - z^2|^{1/2}/|z|$. Then by (5) and the formula above, we get

$$A_{S^{n-1}} + 2\frac{\text{vol}^2(B(0, 1))}{(2\pi)^n} = \frac{2}{\pi n(n-2)!} \int_{\partial \mathbf{K}_+} \frac{|1 - z^2|^{1/2}}{|z|^{n+1}}\, |dz| + \frac{4}{\pi n(n-2)!} \int_1^\infty \frac{\sqrt{t^2 - 1}}{t^{n+1}}\, dt.$$

It remains to show that the second term in the l.h.s. equals the second one in the r.h.s.

After the change of variables $t = 1/s$, and then setting $y_1^2 + \ldots y_{n-1}^2 = s^2$, we get

$$\int_1^\infty \frac{\sqrt{t^2 - 1}}{t^{n+1}}\, dt = \int_0^1 \sqrt{1 - s^2}\, s^{n-2}\, ds = \frac{1}{2\omega_{n-1}} \int_{y_1^2 + \ldots + y_n^2 \leq 1}\, dy = \frac{\omega_n}{2n\omega_{n-1}} = \frac{\sqrt{\pi}\Gamma(\frac{n-1}{2})}{2n\Gamma(\frac{n}{2})},$$

where $\omega_n = 2\pi^{n/2}/\Gamma(n/2)$ is the area of $S^{n-1}$ (and then $\text{vol}(B(0, 1)) = \omega_n/n$). The proof is then reduced to showing that

$$\frac{2}{\pi n(n-2)!} \frac{\sqrt{\pi}\Gamma(\frac{n-1}{2})}{2n\Gamma(\frac{n}{2})} = \frac{\text{vol}^2(B(0, 1))}{(2\pi)^n},$$

and that can be verified by uisng the Legendre duplication formula $\Gamma(2t) = (2\pi)^{-1/2} 2^{2t-1/2} \Gamma(t)\Gamma(t + 1/2)$ applied to $t = (n - 1)/2$. $\square$

# 6 Proof of Theorem 2 and Corollary 1

*Proof of Theorem 2:* The proof of Theorem 2 follows directly from Proposition 2 and Theorem 5(b). We integrate (42) w.r.t. $\theta \in [0, \pi]$, and then pass from polar coordinates to Cartesian ones.

*Proof of Corollary 1:*

Consider (i). By the standard Weyl asymptotics, $N^\sharp(r) = \tau_n(R^n - \text{vol}(\mathcal{O}))r^n(1 + O_R(1/r))$. Apply Theorem 2 and in $\limsup M(r)/r^n$, take the limit $R \to R_0+$.

Next, consider (ii). For $N^\sharp$, we have $N^\sharp(r) = \tau_n R^n r^n(1 + O_R(1/r))$. Take the limit $R \to R_0+$ as above to prove the corollaty in case (ii).

Finally, in case (iii), we have $N^\sharp(r) = (2\pi)^{-n} \int_{|x| \leq R,\, \sum c^2 g^{ij} \xi_i \xi_j \leq 1}\, dx\, d\xi(1 + O_R(1/r))$, and as above, we get the desired estimate.

# 7 Asymptotics of the sphere resonances, proof of Theorem 1

It is enough to study the case $R_0 = 1$. It is well known (and also follows from section 4) that the resonances of the unit sphere $S^{n-1}$, $n$ odd, coincide with the zeros of $H_\nu^{(1)}(\lambda)$, $\nu = l + n/2 - 1$, $l = 0, 1, \ldots$, and each one has multiplicity $m(l)$, see (28). By [O2] and (38), they lie in an $O(1/\nu)$ neighborhood of $\nu\partial \mathbf{K} \cap \{\Im z < 0\}$ and are symmetric to the zeros of $H_\nu^{(2)}(\lambda)$ about the real axis, see Figure 2. They are the zeros of the polynomial $e^{-i\lambda}\lambda^\nu H_\nu^{(2)}(\lambda)$ of degree $\nu - 1/2$. More precisely, one can describe those resonances as follows.

Set

$$a_k = \left[\frac{3}{2}\left(-\frac{\pi}{4} + k\pi\right)\right]^{2/3}, \quad k = 1, 2, \ldots.$$

We view $a_k$ as approximate zeros of $\mathrm{Ai}(-z)$, and by (34), the actual zeros of $\mathrm{Ai}(-w)$ stay at distance no greater than $C/\sqrt{a_k}$. Motivated by (38), we then set

$$\lambda_{\nu k} = \nu \zeta^{-1}(\nu^{-2/3} e^{-i\pi/3} a_k) = \nu \rho^{-1}(-i(k-1/4)\pi/\nu), \quad k = 1, 2, \ldots, \nu - \frac{1}{2}. \tag{57}$$

We view $\bar{\lambda}_{\nu k}$, $k = 1, \ldots, l + (n-3)/2$, $l = 0, 1, \ldots$ (recall (47)), as approximate resonances. The difference between $\lambda_{\nu k}$ and the zeros $h_{\nu k}$ of $H_\nu^{(2)}$ can be estimated as follows. We combine asymptotics (34), (35) of Ai, Ai$'$ with that of $H_\nu^{(2)}$ (38). Then we get

$$\cos(i\nu\rho - \pi/4) + O(1/\nu|\rho|) = 0, \tag{58}$$

where $\rho = \rho(\lambda/\nu)$ and ignoring the remainder, we do get solutions $\lambda_{k\nu}$ as in (57). To simplify our analysis, we will analyze only the zeros $h_{\nu k}$ and the approximate ones $\lambda_{k\nu}$ in the sector

$$\epsilon \leq \arg \lambda \leq \pi - \epsilon, \tag{59}$$

with $0 < \epsilon \ll 1$ fixed, and we will estimate roughly the rest of the zeros. This will guarantee, that we work with $\rho$ that is away from 0 and $-i\pi$. One can, in principle, work in the whole sector $\arg \lambda \in [0, \pi]$ and use the fact that either though $\rho$ is not analytic (and invertible) near 1, the function $\zeta$ is; and this would probably give a sharp bound on the remainder term in Theorem 1.

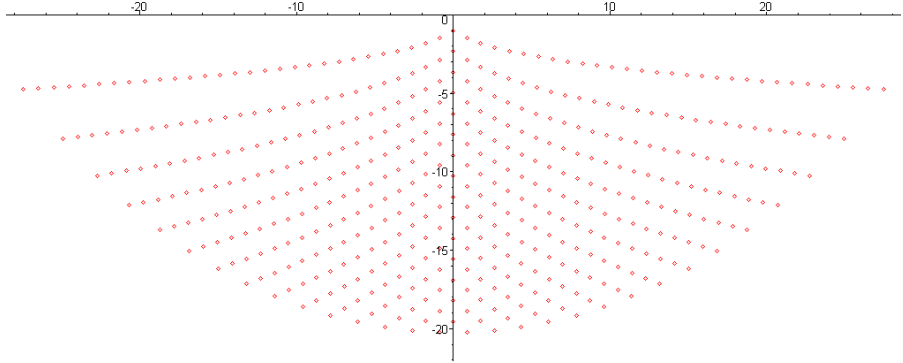

Figure 2: Resonances of the unit sphere $S^2$ (with Dirichlet b. c.), corresponding to $\nu \leq 30$.

For $\lambda$ as in (59), $\rho = \rho(\lambda/\nu)$ satisfies $|\rho| \geq \epsilon^{3/2}/C$. Since the spacing between distinct $\lambda_{\nu k}$ is uniformly bounded from below in (59), by the Rouché theorem we get from the equation $\cos(i\nu\rho(\lambda/\nu) - \pi/4) + O(1/\nu) = 0$ for $h_{\nu k} = \lambda$ that $|\lambda_{\nu k} - h_{\nu k}| \leq C/\nu$ for those $h_{\nu k}$ in the sector (59). For the remaining $h_{\nu k}$ we have

$$\#\{h_{\nu k}; \ \arg h_{\nu k} \leq \epsilon \text{ or } \arg h_{\nu k} \geq \pi - \epsilon\} \leq C\epsilon^{3/2}\nu, \tag{60}$$

and they stay at a distance $O(1/\nu)$ from $\partial \mathbf{K}_+$ in the sector (59).

We will estimate now the number

$$\Lambda(r, \theta, \theta + \Delta\theta) = \sum_{|\lambda_{\nu k}| \leq r; \ \theta \leq \arg \lambda_{\nu k} \leq \theta + \Delta\theta} m(l), \tag{61}$$

see (28), for $\theta$ in (59) and $0 < \Delta\theta$ small enough. Let $z(\theta)$ be defined by $z \in \partial \mathbf{K}$, $\arg z = \theta$. Then $[0, \pi] \ni \theta \mapsto z(\theta)$ is a parameterization of $\partial \mathbf{K}_+$. On the other hand, the properties of $\rho$ imply that $[0, \pi] \ni t \to \rho^{-1}(-it) = z$ is another parameterization. Differentiating $\rho(z) = -it$, we get $\rho'(z)(dz/dt) = -i$, therefore

$$\frac{dt}{dz} = i\rho'(z) = -i\frac{\sqrt{1-z^2}}{z}. \tag{62}$$

Note that for any $\lambda_{\nu k}$ appearing in (61), $|\lambda_{\nu k}| = \nu|z(\theta)| + \nu O(\Delta\theta)$. On the other hand, for a fixed $\nu$, the number $\Delta k$ of $k$'s appearing in (61) satisfies $\pi \Delta k/\nu = \Delta t + O(1/\nu)$, therefore the number of $\lambda_{\nu k}$'s corresponding to the interval $[t, t + \Delta t]$ is $\nu \Delta t + O(1)$. Therefore,

$$\begin{aligned}
\Lambda(r, \theta, \theta + \Delta\theta) &= \sum_{\nu|z(\theta)| \le r} m(l)\frac{\nu \Delta t}{\pi} + r^n O((\Delta t)^2) + O(r^{n-1}) \\
&= \sum_{\nu \le r/|z(\theta)|} \frac{2\nu^{n-1}}{\pi(n-2)!}\Delta t + r^n O((\Delta t)^2) + O(r^{n-1}).
\end{aligned}$$

Note that for $\theta$ as in (59), $t$ stays at distance to the endpoints $0$ and $\pi$ of $[0, \pi]$ at least $\sqrt{\epsilon}/C$. This, in combination with (62), and the inequality above, yields

$$\begin{aligned}
\Lambda(r, \epsilon, \pi - \epsilon) &= \frac{2r^n}{\pi n(n-2)!} \int_{\partial \mathbf{K}_+} \frac{dt}{|z|^n} + r^n O(\epsilon) + O(r^{n-1}) \\
&= \frac{2r^n}{\pi n(n-2)!} \int_{\partial \mathbf{K}_+} \frac{|1-z^2|^{\frac{1}{2}}}{|z|^{n+1}} d|z| + r^n O(\epsilon) + O(r^{n-1}),
\end{aligned}$$

compare with Lemma 4. By (60),

$$\Lambda(r, 0, \epsilon) + \Lambda(r, \pi - \epsilon, \pi) \le C\epsilon^{3/2} r^n.$$

The estimates above, true for any $0 < \epsilon \ll 1$, combined with Lemma 4, prove Theorem 1 for $\Lambda(r, 0, 2\pi)$. The relation between $\lambda_{\nu k}$ and $h_{\nu k}$ (and the resonances $\bar{h}_{\nu k}$) established above yields the same result for $N_{S^{n-1}}(r)$. $\square$

## 8 Proof of Theorem 3

As mentioned in the Introduction, the asymptotic formula $N(r) = K_n R_0^n r^n = o(r^n)$ was proven in [Z2], so we only need to show that $K_n$ is given by (5). As above, it is enough to study the case $R_0 = 1$. Note first that in the notation of [Z2], $\xi$ there corresponds to $-i\rho$ here, therefore, $\Re\xi = \Im\rho$. The approach in [Z2] is to find approximate $\rho_{\nu k}$ in Lemma 6 there:

$$\rho_{\nu k} = i\frac{\pi k}{\nu} + \frac{\log \nu}{\nu} + \frac{1}{\nu} \log f\left(\frac{\pi k}{\nu}\right), \quad k = -[\nu/2], \ldots, -1, 0, 1, \ldots, \tag{63}$$

$f(\rho) = 1 - z^2$, and to consider "approximate resonances" $\lambda_{\mu k}$ given as solutions of $\rho(\lambda/\nu) = \rho_{\nu k}$. For $k < 0$, $\lambda_{\mu k}$ lie near $\nu \partial \mathbf{K}_+ \cap \{\Re z > 0\}$ (at a distance not exceeding $O(\log \nu)$) similarly to the zeros of $h_\nu^{(1)}$, and those with negative real parts are symmetric to them. For $k > 0$, they lie in a logarithmic neighborhood of the positive real axis, and again, there is a sequence with negative real parts, symmetric to them. They can be considered as approximate zeros of $j_\nu(\lambda)$. Each of the $\lambda_{\mu k}$'s is counted with multiplicity $m(l)$. The counting function of the first type is denoted in [Z2] by $n_-(r)$, the one related to $k > 0$ is $n_+(r)$. Then one replaces $\rho_{\nu k}$ in (63) by its first term $i\pi k/\nu$ only and counts the resulting $\lambda'$s. The proof of [Z2] shows that this gives the leading term in the asymptotics for the exact resonances. On the other hand, the counting function $n_-(r)$ for $\rho^{-1}(i\pi k/\nu)$ with $k < 0$ has the same asymptotics as that of $N_{S^{n-1}}(r)$, see (57). Similarly, the counting function $n_-(r)$ for $\rho^{-1}(i\pi k/\nu)$ with $k > 0$ coincides with that of the zeros of $j_\nu(\lambda)$, counted with multiplicities as can be shown as in (58). The latter is governed by the classical Weyl law, as it represents the eigenvalues of the Dirichlet problem in $B(0, 1)$. Therefore,

$$N(r) = N_{S^{n-1}}(r) + \frac{\text{vol}^2(B(0,1))}{(2\pi)^n} r^n + o(r^n),$$

which proves Theorem 8. We also note that one could also compare the integrals at the end of the proof of [Z2] on p. 402, and compare them with the integral representation of $A_{S^{n-1}}$ we have, see also the proof of Lemma 4, and to deduce the proof from there. $\square$

## 9 Proof of Theorem 4

Let us mention first, that by the Weyl formula for the transmission problem,

$$N^\sharp(r) = \frac{1}{(2\pi)^n} \int_{|x| \le R_0;\, c^2|\xi|^2 \le 1} \mathrm{d}x\, \mathrm{d}\xi + \tau_n(R^n - R_0^n)r^n + o(r^n), \tag{64}$$

where $R$ is related to the reference operator $P^\sharp$, as in the Introduction. This shows that the statement of the theorem shows that the estimate in Theorem 2 turns into asymptotic, when $R \to R_0+$.

We will next find an explicit expression of $s(\lambda)$ in this case. Set

$$j_l(t) = t^{1-n/2} J_{l+n/2-1}(t), \quad h_l^{(1,2)}(t) = t^{1-n/2} H_{l+n/2-1}^{(1,2)}(t).$$

Without loss of generality, we may assume that $R_0 = 1$. Let $u \in H^2_{\text{loc}}(\mathbf{R}^n)$ be any solution of $(P - \lambda^2)u = 0$. Write $u$ in polar coordinates, and project onto the eigenspace spanned by the spherical harmonics corresponding to the momentum $l$ to get for the following for the projection $u_l$:

$$u_l(r\omega) = a_l(\lambda) j_l(\lambda r/c), \quad r < 1.$$

Outside $B(0, 1)$, we have

$$u_l(r\omega) = b_l^{(1)}(\lambda) h_l^{(1)}(\lambda r) + b_l^{(2)}(\lambda) h_l^{(2)}(\lambda r), \quad r > 1.$$

By the transmission conditions implied by the requirement $u \in H^2_{\text{loc}}(\mathbf{R}^n)$, we get

$$b_l^{(1)}(\lambda) h_l^{(1)}(\lambda) + b_l^{(2)}(\lambda) h_l^{(2)}(\lambda) = a_l(\lambda) j_l(\lambda/c), \tag{65}$$
$$b_l^{(1)}(\lambda) h_l^{(1)'}(\lambda) + b_l^{(2)}(\lambda) h_l^{(2)'}(\lambda) = c^{-1} a_l(\lambda) j_l'(\lambda/c). \tag{66}$$

Solving the system above, we get

$$\frac{b_l^{(2)}(\lambda)}{b_l^{(1)}(\lambda)} = -\frac{c h_l^{(2)'}(\lambda) j_l(\lambda/c) - h_l^{(2)}(\lambda) j_l'(\lambda/c)}{c h_l^{(1)'}(\lambda) j_l(\lambda/c) - h_l^{(1)}(\lambda) j_l'(\lambda/c)}.$$

The quotient above is the absolute scattering matrix acting on the space spanned by the spherical harmonics with momentum $l$. Since the absolute values of the determinants of the relative and absolute scattering matrices coincide, we see that

$$|s(\lambda)| = \prod_{l=0}^\infty \left| \frac{c h_l^{(2)'}(\lambda) j_l(\lambda/c) - h_l^{(2)}(\lambda) j_l'(\lambda/c)}{c h_l^{(1)'}(\lambda) j_l(\lambda/c) - h_l^{(1)}(\lambda) j_l'(\lambda/c)} \right|^{m(l)}.$$

This characterizes the resonances as the zeros of the denominator above with multiplicities $m(l)$. We will show that they split into two groups: one near the real axis near the zeros of $j_l(\lambda/c)$ (or those of $j_l'(\lambda/c)$ that have the same asymptotic); and another one near the zeros of $h_l^{(1)}(\lambda)$ (we can also view them as approximating the zeros of $h_l^{(1)'}(\lambda)$). The conjugates of those resonances in $\mathbf{C}_+$ coincide with the zeros of the numerator, so we will study those zeros instead of the resonances.

Therefore, we are interested in the asymptotic distribution of the solutions of

$$ch_l^{(2)'}(\lambda) j_l(\lambda/c) - h_l^{(2)}(\lambda) j_l'(\lambda/c) = 0, \tag{67}$$

each one counted $m(l)$ times. This equation is equivalent to

$$\frac{\lambda h_l^{(2)'}(\lambda)}{h_l^{(2)}(\lambda)} - \frac{(\lambda/c) j_l'(\lambda/c)}{j_l(\lambda/c)} = 0, \tag{68}$$

that can be viewed as the transmission condition satisfied by the resonant states at the resonance frequencies. Note that the terms above are regular at $\lambda = 0$. We will use below the notation $\zeta = \zeta(\lambda/\nu)$, $\zeta_c = \zeta(\lambda/c\nu)$, where $\nu$ is as in (47). We also reserve the notation $z$ for $\lambda/\nu$. Fix $0 < \epsilon \ll 1$. We will study the zeros of the equation above for $z = \lambda/\nu \in \Omega_\epsilon := \mathbf{C} \setminus \{D(c, \epsilon) \cup D(1, \epsilon)\}$, (that must be lie $\mathbf{C}_+$), and then roughly estimate the zeors in the two disks appearing above. Here, $D(a, r)$ is the disk in $\mathbf{C}$ centered at $a$ with radius $r$. Note that the exclusion of a neighborhood of $z = 1$ removes a neighborhood of $\zeta = 0$. On the other hand, $z \notin D(c, \epsilon)$ guarantees that $|\zeta_1| \geq 1/C$. We will first express the ratio $h_l^{(2)'}(\lambda)/h_l^{(2)}(\lambda)$ away from its zeros and poles that lie near $\nu \partial \mathbf{K}_+$. To use the asymptotics (38) and (33), arg $\zeta$ must be outside a fixed neighborhood of $-\pi/3$. In $\Omega_\epsilon \cap \mathbf{C}_+$, this is achieved if $z$ is at a fixed distance from $\partial \mathbf{K}_+$. Therefore,

$$\frac{\lambda}{\nu} \frac{h_l^{(2)'}(\lambda)}{h_l^{(2)}(\lambda)} = \frac{\lambda}{\nu} \zeta' \frac{\text{Ai}'_+(\nu^{2/3}\zeta)}{\nu^{1/3} \text{Ai}_+(\nu^{2/3}\zeta)} + O\left(\frac{1}{\nu}\right) \tag{69}$$

$$= \pm\sqrt{1 - z^2} + O\left(\frac{1}{\nu}\right) \quad \text{for } \text{dist}(z, \partial \mathbf{K}_+) \geq \epsilon, |\arg z| \leq \pi - \epsilon, \tag{70}$$

where we choose the positive sign for $z \notin \mathbf{K}$, and the negative one otherwise. The last restriction is not significant because of the symmetry of the resonances about the imaginary axis. Above, $\zeta' = \zeta'(\lambda/\nu)$, and the branch of $\sqrt{1 - z^2}$ is chosen so that its imaginary part is non-positive for $\Im z \geq 0$. Similarly,

$$\frac{\lambda j_l'(\lambda/c)}{\nu c j_l(\lambda/c)} = \frac{\lambda \zeta_c'}{\nu c} \frac{\text{Ai}'(\nu^{2/3}\zeta_c)}{\nu^{1/3} \text{Ai}(\nu^{2/3}\zeta_c)} + O\left(\frac{1}{\nu}\right) \tag{71}$$

$$= \sqrt{1 - z^2/c^2} + O\left(\frac{1}{\nu}\right) \quad \text{for } \epsilon \leq |\arg z| \leq \pi - \epsilon. \tag{72}$$

It is easy to see that the leading terms in (70) and (72) are never equal in the common region of validity, moreover, the absolute value of their difference is bounded from below by a positive constant. Therefore, for $\nu \gg 1$, equation (68) may have solutions only in

$$\left( \{0 \leq \arg z \leq \epsilon\} \cup \{\pi - \epsilon \leq \arg z \leq \pi\} \right) \cup \{\text{dist}(z, \partial \mathbf{K}_+) \leq \epsilon\}. \tag{73}$$

By the symmetry of resonances about the imaginary axis, it is enough to study them in $\Re \lambda \geq 0$ only.

Let us first focus on the region $\Omega_1 = \{\text{dist}(z, \partial \mathbf{K}_+) \leq \epsilon\} \cap \{|z - 1| \geq \epsilon; \Re z > 0\}$. Isolating a neighborhood of $z = 1$ allows us to use the asymptotics of Ai since $|\zeta| > 1/C$. Then (69) is still valid but to get an analogue of (70), we need to use (34), (35) instead of (33). We divide (67) by $j_l(\lambda/c)$ and after canceling an elliptic factor, we write it in the form

$$0 = i\sqrt{1 - z^2} \sin\left(i\nu\rho - \frac{\pi}{4}\right)\left(1 + O\left(\frac{1}{\nu}\right)\right) \tag{74}$$

$$- \left(\sqrt{1 - z^2/c^2} + O\left(\frac{1}{\nu}\right)\right) \cos\left(i\nu\rho - \frac{\pi}{4}\right)\left(1 + O\left(\frac{1}{\nu}\right)\right).$$

This can be written also as
$$\tan\left(i\nu\rho - \frac{\pi}{4}\right) = -i\frac{\sqrt{1 - z^2/c^2}}{\sqrt{1 - z^2}} + O\left(\frac{1}{\nu}\right).$$

For $w \neq \pm i$, the equation $\tan z = w$ has a unique solution $z = \arctan w$ in $0 \leq \Re z < \pi$, and by the periodicity of $\tan z$, all solutions are given by shifts by $k\pi, k = 0, \pm 1, \ldots$. For $\nu \gg 1$, the r.h.s above is at uniformly bounded from below, at positive distance to $\pm i$ for $z \in \Omega_1$, therefore we get the equation
$$i\rho = \frac{1}{\nu}\arctan\frac{\sqrt{1 - z^2/c^2}}{i\sqrt{1 - z^2}} + \frac{k + 1/4}{\nu}\pi + O\left(\frac{1}{\nu^2}\right).$$

We think of this as an equation for $\rho$, where $z = z(\rho)$. Since $z \in \Omega_1$, we wee that $C\epsilon^{3/2}\nu \leq k \leq \pi/2 + O(1/\nu)$, see (60). For $\nu \gg 1$, we get solutions
$$\rho_{\nu k} = -\frac{k\pi}{\nu}i + O\left(\frac{1}{\nu}\right)$$

for $k$ as above. As in section 7, they are mapped into $z$'s approximating zeros of $H_\nu^{(2)}(\nu z)$, with deviation $O(1/\nu)$, and therefore we get conjugate resonances $\lambda_{\nu k}$ in an $O(1)$ neighborhoods of those zeros of $h_\nu^{(2)}(\lambda)$ that lie in an $O(1)$ neighborhood of $\nu\Omega_1$. The $O(1)$ error would not change the asymptotic of their counting function, therefore we get
$$\#\{\lambda_{\nu k}; |\lambda_{\nu k}| < r; \nu \geq \nu_0\} = \frac{1}{2}\left(A_{S^{n-1}} + O(\epsilon^{3/2})\right) r^n(1 + o(1)) \tag{75}$$

based on the analysis in section 7, where $\nu_0 > 0$ is fixed. The factor $1/2$ there is explained by the fact that we work in $\Re\lambda > 0$ only.

We next study the solutions of (68) in $\Omega_2 = \{0 \leq \arg z \leq \epsilon; |z - 1| \geq \epsilon; |z - c| \geq \epsilon\}$. Similarly to the analysis above, the reason to remove a neighborhood of $z = c$ is to be ensure that $|\zeta_c(z)| = |\zeta(z/c)| \geq 1/C$, therefore (71) is still valid. This case is analyzed in a manner similar to that above, we will skip the details. As a result, if $c > 1$, we get zeros in $\Omega_2$ corresponding to
$$\rho_{c,\nu k} = \frac{k\pi}{\nu}i + O\left(\frac{1}{\nu}\right), \quad k \geq C\epsilon^{3/2}\nu, \quad \nu \geq \nu_0, \tag{76}$$

and $k$ is an integer, compare with (63). The corresponding $\tilde{\lambda}_{\nu k} = \nu c \rho^{-1}(\rho_{c,\nu k})$ approximate zeros of $j_l(\lambda)$ and have counting function satisfying
$$\#\{\tilde{\lambda}_{\nu k}; |\tilde{\lambda}_{\nu k}| < r; \nu \geq \nu_0\} = (\tau_n c^{-n} + O(\epsilon)) r^n(1 + o(1)). \tag{77}$$

Note that the corresponding $z_{c,\nu k} = \lambda_{c,\nu k}/\nu$ are in $[c, \infty) + O(1/\nu)$, thus the restriction $|z - 1| \geq \epsilon$ does not play any role in this case ($c > 1$). If $c < 1$, we have to remove $O(\epsilon)\nu$ number of $k$'s from (76), and (77) will still be preserved.

It remains to estimate the number $N_{1,\nu}(\epsilon)$ and $N_{c,\nu}(\epsilon)$ of zeros of (68) in the discs $D(1, \epsilon)$ and $D(c, \epsilon)$ that are not covered by $\Omega_1 \cup \Omega_2$. We will prove that
$$N_{1,\nu}(\epsilon) + N_{c,\nu}(\epsilon) \leq C\epsilon\nu + C_\epsilon \log \nu. \tag{78}$$

Set
$$f_\nu(z) = ch_l^{(2)'}(\nu z) j_l(\nu z/c) - h_l^{(2)}(\nu z) j_l'(\nu z/c).$$

We will estimate the number $N_{1-i\epsilon,\nu}(2\epsilon)$ of zeros in $|z - (1 - i\epsilon)| \leq 2\epsilon$ and then use the fact that $N_{1,\nu} \leq \tilde{N}_{1,\nu}$. By Jensen's formula,
$$\log|f_\nu(1 - i\epsilon)| + \int_0^{4\epsilon} \frac{N_{1-i\epsilon,\nu}(s)}{s} ds = \frac{1}{2\pi}\int_0^{2\pi} \log|f_\nu(1 - i\epsilon + 4\epsilon e^{i\theta})| d\theta.$$

By a standard argument,
$$N_{1-i\epsilon,\nu}(2\epsilon)\log 2 \leq \frac{1}{2\pi}\int_0^{2\pi}\log|f_\nu(1 - i\epsilon + 4\epsilon e^{i\theta})| d\theta - \log|f_\nu(1 - i\epsilon)|. \tag{79}$$

Write $f_\nu(z) = j_l(\nu z/c)h_l^{(2)}(\nu z)g_\nu(z)$, and as in (69), (71), we show that $|g_\nu(1-i\epsilon)| \geq 1/C - C_\epsilon/\nu$. Let $c > 1$. Then by the asymptotics of $J_\nu$, $H_\nu^{(2)}$, $\log|j_l(\nu z/c)h_l^{(2)}(\nu z)| = \nu(-\rho(1/c) + O(\epsilon)) + O(\log \nu)$ on $|z - (1-\epsilon)| = 4\epsilon$, and it is also true for $z = 1 - \epsilon$. Plugging this into (79), we get (78) for $N_{1,\nu}(\epsilon)$. Let $c < 1$. Then $\log|j_l(\nu z/c)h_l^{(2)}(\nu z)| = O(\nu\epsilon) + O(\log \nu)$ on $|z - (1-\epsilon)| = 4\epsilon$, and it is also true for $z = 1 - \epsilon$. Thus we get (78) for $N_{1,\nu}(\epsilon)$ in this case as well. In the same way we prove (78) for $N_{c,\nu}(\epsilon)$.

Estimate (78) show that the zeros missed in (75), (77) affect the leading term in the asymptotic of $N_P(r)$ only by an $O(\epsilon)$ term. On the other hand, summing up (75) and (77), and then taking the limit $\epsilon \to 0$ in $\limsup N(r)/r^n$, completes the proof of Theorem 4. □